\documentclass[11pt,twoside]{article}
\usepackage{amssymb}
\usepackage{amsmath,amssymb}\usepackage{cite}\usepackage{mathrsfs}\usepackage[amsmath,thmmarks]{ntheorem}
\usepackage{listings}
\usepackage[titletoc]{appendix}\allowdisplaybreaks

\makeatletter
\renewcommand\theequation{\thesection.\@arabic\c@equation}
\makeatother

\newtheorem{thm}{Theorem}[section]%
\newtheorem{lem}[thm]{Lemma}%
\newtheorem{Exam}[thm]{Example}%
\newtheorem{cor}[thm]{Corollary}%
\newtheorem{Pro}[thm]{Proposition}%
\newtheorem{Que}[thm]{Question}%
\newenvironment*{proof}{{\bf Proof.}}{\endproof}
\def\endproof{$\Box$}

\input cyracc.def

\begin{document}

\title{{\bf Sign-Balanced Pattern-Avoiding Permutation Classes}}
\footnotetext{$^{*}$ Corresponding author.
\\ E-mail addresses: guopf999@163.com (Pengfei Guo), Junyao$_{-}$Pan@126.com (Junyao Pan).}

\author{{\bf Junyao Pan$^{{\rm 1}}$, Pengfei Guo$^{{\rm 2}, *}$}\\
 \\{\footnotesize $1.$ School of Sciences, Wuxi University, Wuxi, 214105, P. R. China}\\
{\footnotesize $2.$ School of Mathematics and Statistics, Hainan Normal University, Haikou 571158, P. R. China}}

\date{}
\maketitle

\noindent {\bf Abstract:} A set of permutations is called sign-balanced if the set contains the same number of even permutations as odd permutations.
Let $S_n(\sigma_1, \sigma_2, \ldots, \sigma_r)$ be the set of permutations in the symmetric group $S_n$ which avoids patterns $\sigma_1, \sigma_2, \ldots, \sigma_r$.
The aim of this paper is to investigate when, for certain patterns $\sigma_1, \sigma_2, \ldots, \sigma_r$,
$S_n(\sigma_1, \sigma_2, \ldots, \sigma_r)$ is sign-balanced for every integer $n>1$.
We prove that for any $\{\sigma_1, \sigma_2, \ldots, \sigma_r\}\subseteq S_3$,
if $\{\sigma_1, \sigma_2, \ldots, \sigma_r\}$ is sign-balanced except $\{132, 213, 231, 312\}$,
then $S_n(\sigma_1, \sigma_2, \ldots, \sigma_r)$ is sign-balanced for every integer $n>1$.
In addition, we give some results in the case of avoiding some patterns of length $4$.



\section {Introduction}

In this paper, $S_n$ is always the symmetric group of degree $n$, its identity element is used by $id_n$, and the binomial coefficient is denoted by $C^k_n$.

Fix a permutation $\sigma=\sigma_1\sigma_2\cdots \sigma_k\in S_k$ and let $\pi=\pi_1\pi_2\cdots \pi_n\in S_n$ with $k\leq n$.
If there exists a subset of indices $1\leq i_1<i_2<\cdots <i_k\leq n$ such that $\pi_{i_s}>\pi_{i_t}$ if and only if $\sigma_s>\sigma_t$
for any $1\leq s<t\leq k$, then we call that $\pi$ contains \emph{the pattern} $\sigma$, and the subsequence
$\pi_{i_1}\pi_{i_2}\cdots \pi_{i_k}$ is called an \emph{occurrence} of $\sigma$ in $\pi$, and expressed by $\sigma\leq\pi$.
Otherwise, we call that $\pi$ \emph{avoids} $\sigma$.
For instance, $132\leq 24153$ since $253$ is an occurrence of $132$ in $24153$, and $53412$ avoids $132$.
In fact, the investigation of pattern avoidance in permutations began in 1968,
when Knuth \cite{K} introduced a stack-sorting machine and showed that a permutation can be sorted to the identity permutation using this machine if and only if it avoids the pattern $231$.
Henceforth, pattern avoidance in permutations has been studied by many scholars, for some details see \cite{A, Baril, BJ, B, M2, S, V}.

Let $\pi=\pi_1\pi_2\cdots \pi_n$ be a permutation in $S_n$. A pair of indices $i<j$ forms an \emph{inversion} in the permutation $\pi$ if $\pi_i>\pi_j$, otherwise $i<j$ forms a \emph{noninversion} in the permutation $\pi$. Additionally, we denote by $\tau(\pi)$ ($\theta(\pi)$) the number of inversions (noninversions) in $\pi$. If $\tau(\pi)$ is an even number, then we say that $\pi$ is an \emph{even permutation}, and otherwise $\pi$ is an \emph{odd permutation}. In addition, we say that a set of permutations is \emph{sign-balanced} if the number of even permutations equals the number of odd permutations in this set. Let $S_n(\sigma)$ denote the set of permutations in $S_n$ that avoids pattern $\sigma$. Simion and Schmidt \cite{SS} proved that $S_n(321)$ is sign-balanced if $n$ is even, and the number of even permutations in $S_n(321)$ exceeds the number of odd permutations by the Catalan number $C_{\frac{1}{2}(n-1)}$ if $n$ is odd. Afterwards, the sign-balance property of permutation class that avoids one pattern was studied under various conditions in \cite{EFP, EFP1, R, S}, and further the sign-balance of permutation class with respect to certain statistics was studied in \cite{A1, M1}.

In this note, we consider the permutation class that avoids several patterns. Let $S_n(\sigma_1, \sigma_2, \ldots, \sigma_r)$ be the set of permutations in $S_n$ which avoids all patterns $\sigma_1, \sigma_2, \ldots, \sigma_r$. We focus on the following problem and further obtain some results.

\begin{Que}\label{pan1-1}
Are there some patterns $\sigma_1, \sigma_2, \ldots, \sigma_r$ such that $S_n(\sigma_1, \sigma_2, \ldots, \sigma_r)$ is sign-balanced for every integer $n>1$?
\end{Que}

\begin{thm}\label{pan1-2}
Suppose $\sigma_1, \sigma_2, \ldots, \sigma_r\in S_3$. Then, $S_n(\sigma_1, \sigma_2, \ldots , \sigma_r)$ is sign-balanced for every integer $n>1$ if and only if $\{\sigma_1, \sigma_2, \ldots, \sigma_r\}\neq\{132$, $213$, $231$, $312\}$ and $\{\sigma_1, \sigma_2, \ldots, \sigma_r\}$ is sign-balanced.
\end{thm}

\begin{thm}\label{pan1-3}
All of $S_n(1234, 3214)$, $S_n(4321, 4123)$, $S_n(4321, 2341)$, $S_n(1234, 1432)$, $S_n(1243, 2143)$, $S_n(3421, 3412)$, $S_n(4312, 3412)$, $S_n(2134, 2143)$, $S_n(1423, 1432)$,\\ $S_n(3241, 2341)$, $S_n(4132, 4123)$ and $S_n(2314, 3214)$ are sign-balanced for every integer $n>1$.
\end{thm}

\section {Preliminaries}
Recall some notions and notations. Consider $\sigma\in S_l$ and $\pi\in S_m$. The \emph{direct sum} of $\sigma$ and $\pi$ is denote by $\sigma\oplus\pi$, that is,

\[
  \sigma\oplus\pi(i) =
  \begin{cases}
    \sigma(i), &\text{if $1\leq i\leq l$;}\\
	\pi(i)+l, &\text{if $l+1\leq i\leq l+m$.}
  \end{cases}
\]
and the \emph{skew sum} of $\sigma$ and $\pi$ is denote by $\sigma\ominus\pi$, that is,
\[
  \sigma\ominus\pi(i) =
  \begin{cases}
    \sigma(i)+m, &\text{if $1\leq i\leq l$;}\\
	\pi(i-l), &\text{if $l+1\leq i\leq l+m$.}
  \end{cases}
\]
Moreover, for any $\sigma\in S_n$, its reversal $\overline{\sigma}$ is given by $\overline{\sigma}(i)=\sigma(n+1-i)$;
its complement $\sigma^*\in S_n$ is the permutation $\sigma^*(i)=n+1-\sigma(i)$; the inverse $\sigma^{-1}$ is the usual group theoretic inverse permutation.
Let $R=\{\sigma_1, \sigma_2, \ldots, \sigma_r\}$ be a set of permutations in $S_n$. We use $R^{-1}$ and $R^*$ and $\overline{R}$ to denote $\{\sigma^{-1}_1, \sigma^{-1}_2, \ldots, \sigma^{-1}_r\}$ and $\{\sigma^*_1, \sigma^*_2, \ldots, \sigma^*_r\}$ and $\{\overline{\sigma_1}, \overline{\sigma_2}, \ldots,\overline{\sigma_r}\}$ respectively.

Next we give some results which are useful in proving Theorem\ \ref{pan1-2} and Theorem\ \ref{pan1-3} by inductive method.

\begin{lem}\label{pan2-3}
Let $\pi\in S_n$ and $\mu$ be a permutation that is from exchanging the numbers in two different positions of $\pi$. Then the parity of $\pi$ and $\mu$ is opposite.
\end{lem}
\noindent{\bf Proof} It is well-known that the parity of $\tau(\pi)$ and $\tau(\mu)$ is opposite. Thus, the parity of $\pi$ and $\mu$ is opposite. $\hfill{\Box}$

\begin{lem}\label{pan2-5}
Let $\pi=\pi_1\pi_2\cdots \pi_n\in S_n$ and $\mu=\pi_1\pi_2\cdots \pi_{i}(n+1)\pi_{i+1}\cdots \pi_n\in S_{n+1}$ for some $i=0, 1, \ldots , n$. In particular, if $i=0$ then $\mu=(n+1)\pi_1\pi_2\cdots \pi_n\in S_{n+1}$.
If $n-i$ is an even number, then $\pi$ and $\mu$ have the same parity; if $n-i$ is an odd number, then $\pi$ and $\mu$ have the opposite parity.
\end{lem}
\noindent{\bf Proof} Clearly, $\tau(\mu)=\tau(\pi)+n-i$. Thus, if $n-i$ is an even number, then $\pi$ and $\mu$ have the same parity; if $n-i$ is an odd number, then $\pi$ and $\mu$ have the opposite parity. $\hfill{\Box}$

\begin{lem}\label{pan2-7}
Let $R$ be a set of permutations. Then $S_n(R^*)=S^*_n(R)$ and $S_n(\overline{R})=\overline{S_n(R)}$.
\end{lem}
\noindent{\bf Proof} Let $R=\{\sigma_1, \sigma_2, \ldots, \sigma_r\}$ and $\pi\in S^*_n(R)$. If $\pi\notin S_n(R^*)$, then $\pi$ contains $\sigma^*_i$ that is in $R^*$. Thus, $\pi^*$ contains $\sigma_i$, in other words, $\pi^*\notin S_n(R)$. However, $\pi\in S^*_n(R)$ means that $\pi^*\in S_n(R)$, a contradiction. Therefore, $S^*_n(R)\subseteq S_n(R^*)$. Since $R=(R^*)^*$, it follows that $S^*_n(R^*)\subseteq S_n(R)$, and which implies that $S_n(R^*)\subseteq S^*_n(R)$. Hence, $S_n(R^*)=S^*_n(R)$.
Similarly, we can infer that $S_n(\overline{R})=\overline{S_n(R)}$. $\hfill{\Box}$

\begin{lem}\label{pan2-6}
Let $\sigma\in S_n$. If $C^2_n$ is an even number, then $\sigma^*$, $\overline{\sigma}$ and $\sigma$ have the same parity; otherwise, the parity of $\sigma$ is opposite to both of $\sigma^*$ and $\overline{\sigma}$.
\end{lem}
\noindent{\bf Proof} Note that $\theta(\sigma)+\tau(\sigma)=C^2_n$ and $\theta(\sigma)=\tau(\sigma^*)=\tau(\overline{\sigma})$. Thus, if $C^2_n$ is an even number, then $\tau(\sigma^*)$, $\tau(\overline{\sigma})$ and $\tau(\sigma)$ have the same parity, and otherwise the parity of $\tau(\sigma)$ is opposite to both of $\tau(\sigma^*)$ and $\tau(\overline{\sigma})$, as desired. $\hfill{\Box}$

\begin{cor}\label{pan2-55}
Let $R$ be a set of permutations such that $S_n(R)$ is sign-balanced. Then $S_n(R^*)$ and $S_n(\overline{R})$ are sign-balanced.
\end{cor}
\noindent{\bf Proof} Let $\Omega_+=\{\pi\in S_n(R):$ $\pi$ is an even permutation$\}$ and $\Omega_-=\{\pi\in S_n(R):$ $\pi$ is an odd permutation$\}$. Since $S_n(R)$ is sign-balanced, we deduce that $|\Omega_+|=|\Omega_-|$. Applying Lemma \ref{pan2-7}, it follows that $S_n(R^*)=\Omega^*_+\cup\Omega^*_-$. Set $\Delta_+=\{\pi\in S_n(R):$ $\pi$ is an even permutation$\}$ and $\Delta_-=\{\pi\in S_n(R):$ $\pi$ is an odd permutation$\}$. Then by Lemma \ref{pan2-6}, we see that $\Delta_+=\Omega^*_+$ and $\Delta_-=\Omega^*_-$ if $C^2_n$ is an even number, and $\Delta_-=\Omega^*_+$ and $\Delta_+=\Omega^*_-$ if $C^2_n$ is an odd number. In either case, $|\Delta_+|=|\Delta_-|$ holds. Similarly, we can deduce that $S_n(\overline{R})$ is also sign-balanced. $\hfill{\Box}$

\begin{lem}\label{pan2-66}
Let $\sigma\in S_l$ and $\pi\in S_m$. If $ml$ is an even number, then $\tau(\sigma\ominus\pi)$ and $\tau(\sigma)+\tau(\pi)$ have the same parity; if $ml$ is an odd number,
then the parity of $\tau(\sigma\ominus\pi)$ is opposite to that of $\tau(\sigma)+\tau(\pi)$.
\end{lem}
\noindent{\bf Proof} The lemma follows from the fact that $\tau(\sigma\ominus\pi)=\tau(\sigma)+\tau(\pi)+ml$. $\hfill{\Box}$

\begin{lem}\label{pan2-666}\normalfont \cite[Proposition~17]{SS}
If $R \subseteq S_3$ and $|R| = 4$, then\\
(a) $|S_n(R)| = 0$ if $R \supset \{123, 321\}$ and $n \geq 5$;\\
(b) $|S_n(R)| = 2$ if $R \not\supset \{123, 321\}$, $n \geq 2$.\\
More precisely, the permutations counted in (b) are the appropriate two, depending on $R$,
from among the identity, $23\cdot\cdot\cdot n1$, $n(n -1)\cdot\cdot\cdot312$, their reversals and complements. For the values of $n$ omitted in part (a), we have: $S_1(R) = 1$, $S_2(R) = 2$,
$S_4(R)=0$ except for $S_4(123, 321,132,213) = S_4(123, 321, 231, 312) = 1$.
\end{lem}

\begin{lem}\label{pan2-6666}\normalfont \cite[Erd\H{o}s-Szekeres Theorem]{C}
Any sequence of $mk+1$ distinct real numbers contains either an increasing subsequence of $m+1$ terms or a decreasing subsequence of $k+1$ terms.
\end{lem}

\section {Main results}

We first consider Question \ref{pan1-1} under the assumption of $\sigma_1, \sigma_2, \ldots, \sigma_r \in S_3$. In this case, we see that if $S_n(\sigma_1, \sigma_2, \ldots, \sigma_r)$ is sign-balanced for every integer $n>1$, then $\{\sigma_1,\sigma_2,\ldots,\sigma_r\}$ is sign-balanced because $S_3(\sigma_1,\sigma_2,\ldots,\sigma_r)=S_3\setminus\{\sigma_1,\sigma_2,\ldots,\sigma_r\}$. Thus for Theorem\ \ref{pan1-2}, it suffices to check that $S_n(\sigma_1, \sigma_2, \ldots, \sigma_r)$ is sign-balanced for every integer $n>1$ if $\{\sigma_1, \sigma_2, \ldots, \sigma_r\}\neq\{132, 213, 231, 312\}$ and $\{\sigma_1, \sigma_2, \ldots, \sigma_r\}$ is sign-balanced, and further $S_n(132, 213, 231, 312)$ is not sign-balanced for some integer $n>1$. Now we start from the case that $r=2$.

\begin{Pro}\label{pan3-2}
All of $S_n(132, 231)$, $S_n(312, 213)$, $S_n(132, 312)$ and $S_n(231, 213)$ are sign-balanced for every integer $n>1$.
\end{Pro}
\noindent{\bf Proof} Since $132$ is an odd permutation while $231$ is an even permutation,
we infer that $S_n(132, 231)$ is sign-balanced for $n=2, 3$. Next we prove this proposition by induction on $n$. Assume that $S_n(132, 231)$ is sign-balanced for $n=k>1$. Consider $\pi=\pi_1\pi_2\cdots \pi_{k}\pi_{k+1}\in S_{k+1}$. Note that if $\pi\in S_{k+1}(132, 231)$, then $\pi_1=k+1$ or $\pi_{k+1}=k+1$, otherwise $\pi$ contains either $132$ or $231$. Suppose $\pi_1=k+1$. It is clear that $\pi\in S_{k+1}(132, 231)$ if and only if $\pi'=\pi_2\pi_3\cdots \pi_{k+1}\in S_{k}(132, 231)$. Also, if $\pi_{k+1}=k+1$, then $\pi\in S_{k+1}(132, 231)$ if and only if $\pi''=\pi_1\pi_2\cdots \pi_k\in S_{k}(132, 231)$. Therefore, we deduce that $$S_{k+1}(132, 231)=\{1\ominus\mu|\mu\in S_{k}(132, 231)\}\cup\{\mu\oplus1|\mu\in S_{k}(132, 231)\}.$$
Let $\Omega_{+}=\{1\ominus\mu|\mu\in S_{k}(132, 231),\mu$ is an even permutation$\}$ and $\Omega_{-}=\{1\ominus\mu|\mu\in S_{k}(132, 231),\mu$ is an odd permutation$\}$. Then by Lemma \ref{pan2-5} we infer that if $k$ is even, $\Omega_{+}$ is an even permutation set and $\Omega_{-}$ is an odd permutation set; and if $k$ is odd, $\Omega_{+}$ is an odd permutation set and $\Omega_{-}$ is an even permutation set. Applying inductive hypothesis, we deduce that $\{1\ominus\mu|\mu\in S_{k}(132, 231)\}$ is sign-balanced no matter what the parity of $k$ is. Similarly, we see that $\{\mu\oplus1|\mu\in S_{k}(132, 231)\}$ is sign-balanced. Therefore, $S_{k+1}(132, 231)$ is sign-balanced, and so $S_{n}(132, 231)$ is sign-balanced for every integer $n>1$.

Note $\{312, 213\}=\{132^*, 231^*\}$. It follows from Corollary \ref{pan2-55} that $S_n(312, 213)$ is sign-balanced for every integer $n>1$. Additionally, we see that $\{132, 312\}=\{132^{-1}, 231^{-1}\}$
and $\{231, 213\}=\{312^{-1}, 213^{-1}\}$. Using a result in \cite[Lemma~1]{SS}, we deduce that $S_n(132, 312)=S^{-1}_{n}(132, 231)$ and $S_n(231, 213)=S^{-1}_n(312, 213)$. Since the permutations $\sigma$ and $\sigma^{-1}$ have the same parity, it follows that $S_n(132, 312)$ and $S_n(231, 213)$ are sign-balanced for every integer $n>1$. $\hfill{\Box}$

\begin{Pro}\label{pan3-4}
All of $S_n(123, 213)$, $S_n(321, 231)$, $S_n(123, 132)$ and $S_n(321, 312)$ are sign-balanced for every integer $n>1$.
\end{Pro}
\noindent{\bf Proof} It is clear that $S_n(123, 213)$ is sign-balanced for $n=2, 3$ because $123$ is an even permutation while $213$ is an odd permutation.
By induction on $n$, we assume that $S_n(123, 213)$ is sign-balanced for $n=k>1$. Consider $\pi=\pi_1\pi_2\cdots\pi_{k}\pi_{k+1}\in S_{k+1}$. We see that if $\pi\in S_{k+1}(123, 213)$ then $\pi_1=k+1$ or $\pi_2=k+1$, otherwise $\pi$ contains either $123$ or $213$. In addition, it is straightforward to show that if $\pi_1=k+1$ then $\pi\in S_{k+1}(123, 213)$ if and only if $\pi'=\pi_2\pi_3\cdots \pi_{k+1}\in S_{k}(123, 213)$, and if $\pi_2=k+1$ then $\pi\in S_{k+1}(132, 231)$ if and only if $\pi''=\pi_1\pi_3\cdots \pi_{k+1}\in S_{k}(123, 213)$. Proceeding as in the proof of Proposition\ \ref{pan3-2}, we deduce that $S_{k+1}(123, 213)$ is sign-balanced from Lemma \ref{pan2-5} and the inductive assumption. Therefore, $S_{n}(123, 213)$ is sign-balanced for every integer $n>1$.

In addition, we see that $\{321, 231\}=\{123^*, 213^*\}$, $\{123, 132\}=\{321^{-1}, 231^{-1}\}$ and $\{321, 312\}=\{123^{-1}, 213^{-1}\}$. An argument similar to the one used in the proof of Proposition\ \ref{pan3-2} shows that all of $S_n(321, 231)$, $S_n(123, 132)$ and $S_n(321, 312)$ are sign-balanced for every integer $n>1$.
$\hfill{\Box}$

\begin{Pro}\label{pan3-5}
$S_n(123,321)$ is sign-balanced for every integer $n>1$.
\end{Pro}
\noindent{\bf Proof} It is obvious that $321$ is an odd permutation while $123$ is an even permutation, and thus $S_n(123, 321)$ is sign-balanced for $n=2, 3$.
Additionally, it follows from Lemma\ \ref{pan2-6666} that it suffices to consider $S_4(123, 321)$. One easily checks that $S_4(123, 321)\cap A_4=\{2143, 3412\}$ and $S_4(123, 321)\cap(S_4\setminus A_4)=\{3142, 2413\}$, where $A_4$ is the alternating group of degree $4$. The completes the proof. $\hfill{\Box}$

So far, we have verified all situations of $S_n(\sigma_1, \sigma_2)$ under the assumption of $\sigma_1, \sigma_2\in S_3$.
Secondly, we consider $S_n(\sigma_1, \sigma_2, \sigma_3, \sigma_4)$ in case when  $\sigma_1, \sigma_2, \sigma_3, \sigma_4\in S_3$, and obtain the following result.

\begin{Pro}\label{pan3-6}
Let $\{\sigma_1, \sigma_2, \sigma_3, \sigma_4\}\subseteq S_3$. Then, $S_n(\sigma_1, \sigma_2, \sigma_3, \sigma_4)$ is sign-balanced for every integer $n>1$ if and only if $\{\sigma_1, \sigma_2, \sigma_3, \sigma_4\}\neq\{132,213,231,312\}$ and $\{\sigma_1, \sigma_2, \sigma_3, \sigma_4\}$ is sign-balanced.
\end{Pro}
\noindent{\bf Proof}\ If $\{\sigma_1, \sigma_2, \sigma_3, \sigma_4\}=\{132, 213, 231, 312\}$, then $\{\sigma_1,\sigma_2,\sigma_3,\sigma_4\}\not\supset\{123,321\}$ and further $\{\sigma_1,\sigma_2,\sigma_3,\sigma_4\}\cap\{123, 321\}=\emptyset$. Applying Lemma\ \ref{pan2-666}, we deduce that $S_{n}(132, 213, 231, 312)=\{12\cdots n, n\cdots 21\}$. However, Lemma\ \ref{pan2-6} shows that $12\cdots  n$ and $n\cdots 21$ have the same parity if $C^2_n$ is an even number, and so $S_n(132,213,231,312)$ is not sign-balanced when $C^2_n$ is an even number.

Assume that $\{\sigma_1, \sigma_2, \sigma_3, \sigma_4\}\cap\{123, 321\}\neq\emptyset$ and $\{\sigma_1, \sigma_2, \sigma_3, \sigma_4\}$ is sign-balanced. It suffices to prove that $S_n(\sigma_1, \sigma_2, \sigma_3, \sigma_4)$ is sign-balanced for every integer $n>1$. It is obvious that $S_n(\sigma_1, \sigma_2,\sigma_3,\sigma_4)$ is sign-balanced when $n=2, 3$. Next we consider $S_n(\sigma_1, \sigma_2, \sigma_3, \sigma_4)$ for $n>3$.

Consider $|\{\sigma_1, \sigma_2, \sigma_3, \sigma_4\}\cap\{123, 321\}|=2$. Namely $\{\sigma_1, \sigma_2, \sigma_3, \sigma_4\}\supset\{123, 321\}$. According to Lemma\ \ref{pan2-666}, we deduce that $|S_n(\sigma_1, \sigma_2, \sigma_3, \sigma_4)|=0$ for every integer $n\geq 5$. In addition, Lemma\ \ref{pan2-666} shows that $|S_4(\sigma_1, \sigma_2, \sigma_3, \sigma_4)|\neq 0$ if and only if $\{\sigma_1,\sigma_2,\sigma_3,\sigma_4\}=\{123,321,132,213\}$ or $\{123, 321, 231, 312\}$. However, we see that $\{123, 321, 132, 213\}$ and $\{123, 321, 231, 312\}$ are not sign-balanced. Thus, we deduce that $|S_4(\sigma_1, \sigma_2, \sigma_3, \sigma_4)|=0$. Therefore, $S_n(\sigma_1,\sigma_2,\sigma_3,\sigma_4)$ is sign-balanced for every integer $n>1$.

Consider the case that $|\{\sigma_1, \sigma_2, \sigma_3, \sigma_4\}\cap\{123, 321\}|=1$. Note $321=123^*$. By Corollary \ref{pan2-55}, we can assume that $123\notin\{\sigma_1, \sigma_2, \sigma_3, \sigma_4\}$ and $321\in\{\sigma_1, \sigma_2, \sigma_3, \sigma_4\}$. Since $\{\sigma_1, \sigma_2, \sigma_3, \sigma_4\}$ is sign-balanced, it follows that either $123, 132\notin\{\sigma_1, \sigma_2, \sigma_3, \sigma_4\}$ or $123, 213\notin\{\sigma_1, \sigma_2, \sigma_3, \sigma_4\}$. If $123, 132\notin\{\sigma_1, \sigma_2, \sigma_3, \sigma_4\}$,
then $S_n(\sigma_1, \sigma_2,\sigma_3, \sigma_4)=\{1234\cdots  n, 1234\cdots (n-2)n(n-1)\}$ by Lemma\ \ref{pan2-666}. Clearly, $S_n(\sigma_1,\sigma_2,\sigma_3,\sigma_4)$ is sign-balanced in this case. If $123, 213\notin\{\sigma_1, \sigma_2, \sigma_3, \sigma_4\}$, then $S_n(\sigma_1, \sigma_2, \sigma_3, \sigma_4)=\{1234\cdots  n, 2134\cdots  n\}$ by Lemma\ \ref{pan2-666}. In this case, $S_n(\sigma_1, \sigma_2, \sigma_3, \sigma_4)$ is also sign-balanced. The proof of this proposition is completed. $\hfill{\Box}$

\vskip0.4cm
\noindent{\bf Proof of Theorem \ref{pan1-2}} Consider $S_n(\sigma_1, \sigma_2, \sigma_3, \sigma_4, \sigma_5, \sigma_6)=S_n(S_3)$.
It is obvious that $S_n(S_3)$ is sign-balanced for every integer $n>1$ because $S_n(S_3)=\emptyset$ when $n\geq 3$.
So we have solved Question \ref{pan1-1} under the assumption of $\sigma_1, \sigma_2, \ldots, \sigma_r\in S_3$.
At the same time, we derive Theorem \ref{pan1-2} immediately by Propositions \ref{pan3-2}-\ref{pan3-6}.
$\hfill{\Box}$

\vskip0.4cm
Thirdly, we solve Question \ref{pan1-1} in case when $\sigma_1, \sigma_2, \ldots, \sigma_r\in S_4$, and obtain some results as follows.

\begin{Pro}\label{pan3-7}
Let $R=\{\{1234,3214\}, \{4321,4123\}, \{4321,2341\}, \{1234,1432\}\}$. Then for any $\{\sigma_1, \sigma_2\}\in R$, $S_n(\sigma_1, \sigma_2)$ is sign-balanced for every integer $n>1$.
\end{Pro}
\noindent{\bf Proof} Note that $\{4321, 4123\}=\{\overline{1234}, \overline{3214}\}$, $\{4321, 2341\}=$ $\{1234^{*}, 3214^{*}\}$ and $\{1234$, $1432\}=$ $\{\overline{4321}, \overline{2341}\}$. Then by Corollary\ \ref{pan2-55}, it suffices to prove that $S_n(1234, 3214)$ is sign-balanced for every integer $n>1$. Since $1234$ is an even permutation while $3214$ is an odd permutation, it follows that $S_n(1234, 3214)$ is sign-balanced for $n=2, 3, 4$. Assume that $S_n(1234, 3214)$ is sign-balanced for $1<n\leq k$ with $k>3$. Let $X^{(i)}=\{\pi\in {S_{k+1}(1234, 3214)}:\pi(i)=k+1\}$. It is clear that $S_{k+1}(1234, 3214)$ is the disjoint union $\bigcup^{k+1}_{i=1}X^{(i)}$. It follows from Lemma\ \ref{pan2-6666} that $X^{(i)}=\emptyset$ if $i>5$. Therefore, we deduce that $$S_{k+1}(1234, 3214)=\bigcup^{5}_{i=1}X^{(i)}.$$Let $\pi=\pi_1\pi_2\cdots \pi_{k+1}\in S_{k+1}$ with $\pi_i=k+1$ and $1\leq i\leq 3$. Clearly, $\pi\in X^{(i)}$
if and only if $\pi'=\pi_1\pi_2\cdots \pi_{i-1}\pi_{i+1}\cdots \pi_{k+1}\in S_{k}(1234, 3214)$. In particular, if $i=1$ then $\pi'=\pi_2\pi_3\cdots \pi_{k+1}$. Proceeding as in the proof of Proposition\ \ref{pan3-2}, we deduce that all of $X^{(1)}, X^{(2)}$ and $X^{(3)}$ are sign-balanced by Lemma\ \ref{pan2-5} and the inductive assumption. Next we prove that $X^{(4)}$ and $X^{(5)}$ are also sign-balanced.

Consider $\pi=\pi_1\pi_2\cdots \pi_{k+1}\in S_{k+1}$ with $\pi_4=k+1$. Let $\pi'=\pi_1\pi_2\pi_3\pi_5\cdots \pi_{k+1}$. Clearly, if $\pi\in X^{(4)}$ then $\pi'\in S_{k}(1234,3214)$. For convenience, we set
$$A^{(4)}=\{\pi=\pi_1\pi_2\cdots \pi_{k+1}\in S_{k+1}:\pi_4=k+1,\pi_1\pi_2\pi_{3}\pi_{5}\cdots \pi_{k+1}\in S_{k}(1234,3214)\}.$$
Proceeding as in the proof of Proposition\ \ref{pan3-2}, we deduce that $A^{(4)}$ is sign-balanced by Lemma\ \ref{pan2-5} and the inductive assumption. Suppose that $\pi\notin X^{(4)}$ and $\pi'\in S_{k}(1234,3214)$. It is simple to see that $\pi_1\pi_2\pi_3\pi_4$ is an occurrence of $1234$ or $3214$. Moreover, if $\pi_1\pi_2\pi_3\pi_4$ is an occurrence of $1234$, then $\pi_1<\pi_2<\pi_{3}=k$ otherwise $\pi'$ contains $1234$; if $\pi_1\pi_2\pi_3\pi_4$ is an occurrence of $3214$, then $\pi_1=k>\pi_2>\pi_{3}$ otherwise $\pi'$ contains an occurrence of $3214$. Indeed, if $\pi'\in S_{k}(1234,3214)$, then $\pi\notin X^{(4)}$ if and only if $\pi_1<\pi_2<\pi_{3}=k$ or $\pi_1=k>\pi_2>\pi_{3}$. Let $B^{(4)}=\{\pi=\pi_1\pi_2\cdots \pi_{k+1}\in A^{(4)}:\pi_1=k>\pi_2>\pi_{3}\}$ and $C^{(4)}=\{\pi=\pi_1\pi_2\cdots \pi_{k+1}\in A^{(4)} :\pi_1<\pi_2<\pi_{3}=k\}$. Thus, we have $$X^{(4)}=A^{(4)}\setminus (B^{(4)}\cup C^{(4)}).$$ Consider $B^{(4)}\cup C^{(4)}$. We observe that $B^{(4)}=\big\{\pi=\pi_1\pi_2\cdots \pi_{k+1}\in S_{k+1}:\pi_1=k,\pi_2>\pi_{3},\pi_4=k+1,\pi_2\pi_{3}\pi_{5}\cdots \pi_{k+1}\in S_{k-1}(1234,3214)\big\}$ and $C^{(4)}=\big\{\pi=\pi_1\pi_2\cdots \pi_{k+1}:\pi_1<\pi_2,\pi_{3}=k,\pi_4=k+1,\pi_1\pi_{2}\pi_{5}\cdots \pi_{k+1}\in S_{k-1}(1234,3214)\big\}$. Since $k-1$ and $k-3$ have the same parity, it follows that $B^{(4)}\cup C^{(4)}$ is sign-balanced from Lemma\ \ref{pan2-5} and the inductive assumption. Therefore, $X^{(4)}$ is sign-balanced.

Suppose that $\pi=\pi_1\pi_2\cdots \pi_{k+1}\in X^{(5)}$ with $\pi_5=k+1$. Clearly, $\pi_1\pi_2\pi_{3}\pi_{4}$ avoids $123$ and $321$. Therefore, we deduce that either $\pi_1>\pi_2,\pi_{3}>\pi_{4},\pi_{3}>\pi_{1},\pi_{4}>\pi_{2}$ or $\pi_1<\pi_2,\pi_{3}<\pi_{4},\pi_{3}<\pi_{1},\pi_{4}<\pi_{2}$, otherwise $\pi_1\pi_2\pi_{3}\pi_{4}$ contains an occurrence of $123$ or $321$. Moreover, for any $\pi=\pi_1\pi_2\pi_{3}\pi_{4}\pi_{5}\cdots \pi_{k+1}\in X^{(5)}$, we observe that $\pi_1\pi_{3}\pi_2\pi_{4}\pi_{5}\cdots \pi_{k+1}\in X^{(5)}$. In other words, exchanging the entries $2$ and $3$ is a bijection from $X^{(5)}$ to $X^{(5)}$ and further the parity of the image is opposite to that of the original image under this bijection. Thus, $X^{(5)}$ is sign-balanced. According to the above arguments, it follows that $S_{k+1}(1234, 3214)$ is sign-balanced. The proof of this proposition is completed. $\hfill{\Box}$

\begin{Pro}\label{pan3-17}
Let $R=\{\{1243,2143\},\{3421,3412\},\{4312,3412\},\{2134,2143\}\}$. Then for any $\{\sigma_1,\sigma_2\}\in R$, $S_n(\sigma_1,\sigma_2)$ is sign-balanced for every integer $n>1$.
\end{Pro}
\noindent{\bf Proof} Note that $\{3421,3412\}=\{\overline{1243},\overline{2143}\}$, $\{4312,3412\}=\{1243^{*},2143^{*}\}$ and $\{2134,2143\}=\{3421^{*},3412^{*}\}$.
According to Corollary\ \ref{pan2-55}, it suffices to confirm that $S_n(1243,2143)$ is sign-balanced for every integer $n>1$. Since $2143$ is an even permutation while $1243$ is an odd permutation, we have that $S_n(1234,3214)$ is sign-balanced for $n=2,3,4$. By induction on $n$, we assume that $S_n(1243,2143)$ is sign-balanced for $1<n\leq k$ with $k>3$. Now we consider $S_{k+1}(1243,2143)$. Let $X^{(i)}=\{\pi\in {S_{k+1}(1243,2143)}:\pi(i)=k+1\}$. It is easy to see that $S_{k+1}(1243,2143)$ is the disjoint union $\bigcup^{k+1}_{i=1}X^{(i)}$. We claim that $X^{(i)}$ is sign-balanced for $i=1,2,...,k+1$.

Consider $\pi=\pi_1\pi_2\cdots \pi_{k+1}\in S_{k+1}$ with $\pi_1=k+1$. Note that $\pi\in X^{(1)}$ if and only if $\pi_2\cdots \pi_k\in S_{k}(1243,2143)$. By Lemma\ \ref{pan2-5} and the inductive assumption, we obtain that $X^{(1)}$ is sign-balanced. Similarly, we can infer that $X^{(2)}$ and $X^{(k+1)}$ are sign-balanced.

Let $\pi=\pi_1\pi_2\cdots \pi_{k+1}\in X^{(i)}$ with $2<i<k+1$. Then there at most one number in $\{\pi_1,\pi_2,\ldots,\pi_{i-1}\}$ which is smaller than the biggest number in $\{\pi_{i+1},\ldots,\pi_{k+1}\}$, otherwise $\pi$ contains an occurrence of $1243$ or $2143$. Thus, $\{\pi_1,\pi_2,..., \pi_{i-1}\}$ is either $\{k-i+2,k-i+3,...,k\}$ or $\{m,k-i+3,...,k\}$ for some $m\in\{1,2,...,k-i+1\}$. In the case that $\{\pi_1,\pi_2,..., \pi_{i-1}\}=\{k-i+2,k-i+3,...,k\}$, set $$\Omega_1=\big\{\pi=\pi_1\pi_2\cdots \pi_{k+1}\in X^{(i)}:\{\pi_1,\pi_2,..., \pi_{i-1}\}=\{k-i+2,k-i+3,...,k\}\big\}.$$ We see that $\Omega_1=\{(\sigma\oplus1)\ominus\rho:\sigma\in S_{i-1}(1243,2143),\rho\in S_{k-i+1}(1243,2143)\}$. Proceeding as in the proof of Proposition\ \ref{pan3-2}, we deduce that $\Omega_1$ is sign-balanced by Lemma \ref{pan2-5}, Lemma \ref{pan2-66} and the inductive assumption. In the case that $\{\pi_1,\pi_2,..., \pi_{i-1}\}=\{m,k-i+3,...,k\}$ for every $m\in\{1,2,...,k-i+1\}$, we set $$\Omega^m=\big\{\pi=\pi_1\pi_2\cdots \pi_{k+1}\in X^{(i)}:\{\pi_1,\pi_2,..., \pi_{i-1}\}=\{m,k-i+3,...,k\}\big\}.$$ Notice that exchanging the entries that are $k-i+2$ and $m$ is a bijection from $\Omega^m$ to $\Omega_1$, and further the parity of the image is opposite to that of the original image under this bijection. Since $\Omega_1$ is sign-balanced, it follows that $\Omega^m$ is sign-balanced. Therefore, $X^{(i)}$ is sign-balanced by the fact that $X^{(i)}$ is the disjoint union $\Omega_1\bigcup(\bigcup^{k-i+1}_{m=1}\Omega^m)$. According to the above arguments, it follows that $S_{k+1}(1243,2143)$ is sign-balanced. The proof of this proposition is completed.
$\hfill{\Box}$

\begin{Pro}\label{pan3-8}
Let $R=\{\{1423,1432\},\{3241,2341\},\{4132,4123\},\{2314,3214\}\}$. Then for any $\{\sigma_1,\sigma_2\}\in R$, $S_n(\sigma_1,\sigma_2)$ is sign-balanced for every integer $n>1$.
\end{Pro}
\noindent{\bf Proof} Note that $\{3241,2341\}=\{\overline{1423},\overline{1432}\}$, $\{4132,4123\}=\{1423^*,1432^*\}$ and $\{2314,3214\}=\{3241^{*},2341^{*}\}$.
So it suffices to prove that $S_n(1423,1432)$ is sign-balanced for every integer $1< n\leq4$. Since $1423$ is an even permutation while $1432$ is an odd permutation, we have that $S_n(1423,1432)$ is sign-balanced for $n=2,3,4$. By induction on $n$, we assume that $S_n(1423,1432)$ is sign-balanced for $1<n\leq k$ with $k>3$. Now we consider $S_{k+1}(1423,1432)$. Let $X^{(i)}=\{\pi\in {S_{k+1}(1423,1432)}:\pi(i)=k+1\}$. It is easy to see that $S_{k+1}(1423,1432)$ is the disjoint union $\bigcup^{k+1}_{i=1}X^{(i)}$. We claim that $X^{(i)}$ is sign-balanced for $i=1,2,...,k+1$.

Consider $\pi=\pi_1\pi_2\cdots \pi_{k+1}\in S_{k+1}$ with $\pi_1=k+1$. Note that $\pi\in X^{(1)}$ if and only if $\pi_2\cdots \pi_k\in S_{k}(1423,1432)$. Applying Corollary\ \ref{pan2-5} and the inductive assumption, we obtain that $X^{(1)}$ is sign-balanced. Similarly, we can infer that $X^{(k)}$ and $X^{(k+1)}$ are sign-balanced.

Let $\pi=\pi_1\pi_2\cdots \pi_{k+1}\in X^{(i)}$ with $1<i<k$. Then there at most one number in $\{\pi_{i+1},\pi_{i+2},\ldots,\pi_{k+1}\}$ which is bigger than the smallest number in $\{\pi_{1},\ldots,\pi_{i-1}\}$, otherwise $\pi$ contains an occurrence of $1423$ or $1432$. Therefore, $\{\pi_{i+1},\pi_{i+2},\ldots,\pi_{k+1}\}$ is either $\{1,2....,k-i+1\}$ or $\{1,2....,k-i,m\}$ for some $m\in\{k-i+2,k-i+3,...,k\}$. In the case that $\{\pi_{i+1},\pi_{i+2},\ldots,\pi_{k+1}\}=\{1,2....,k-i+1\}$, set $$\Omega_1=\big\{\pi=\pi_1\pi_2\cdots \pi_{k+1}\in X^{(i)}:\{\pi_{i+1},\pi_{i+2},\ldots,\pi_{k+1}\}=\{1,2....,k-i+1\}\big\}.$$ We see that $\Omega_1=\{(\sigma\oplus1)\ominus\rho:\sigma\in S_{i-1}(1243,2143),\rho\in S_{k-i+1}(1243,2143)\}$. Proceeding as in the proof of Proposition\ \ref{pan3-2}, we deduce that $\Omega_1$ is sign-balanced by Lemma\ \ref{pan2-5} and Lemma\ \ref{pan2-66} and the inductive assumption. In the case that $\{\pi_{i+1},\pi_{i+2},\ldots,\pi_{k+1}\}=\{1,2....,k-i,m\}$ for every $m\in\{k-i+2,k-i+3,...,k\}$, we set $$\Omega^m=\big\{\pi=\pi_1\pi_2\cdots \pi_{k+1}\in X^{(i)}:\{\pi_{i+1},\pi_{i+2},\ldots,\pi_{k+1}\}=\{1,2....,k-i,m\}\big\}.$$ Notice that exchanging the entries that are $k-i+1$ and $m$ is a bijection from $\Omega^m$ to $\Omega_1$, and further the parity of the image is opposite to that of the original image under this bijection. Since $\Omega_1$ is sign-balanced, it follows that $\Omega^m$ is sign-balanced. Therefore, $X^{(i)}$ is sign-balanced by the fact that $X^{(i)}$ is the disjoint union $\Omega_1\bigcup(\bigcup^{k-i+1}_{m=1}\Omega^m)$. According to the above arguments, it follows that $S_{k+1}(1243,2143)$ is sign-balanced. The proof of this proposition is completed.
$\hfill{\Box}$

\vskip0.4cm
\noindent{\bf Proof of Theorem \ref{pan1-3}} Applying Propositions \ref{pan3-7}-\ref{pan3-8}, we have that Theorem \ref{pan1-3} holds immediately.
$\hfill{\Box}$

\vskip0.4cm
However, we observe that $S_n(\sigma_1,\sigma_2)$ is not always sign-balanced for $\sigma_1,\sigma_2\in S_4$.
Finally, we give an example which seems to mean that Question\ \ref{pan1-1} is more and more complicated as the length of avoiding pattern increases.

\begin{Exam}\label{pan3-9}
Both $S_n(1324,2143)$ and $S_n(4231,3412)$ are not sign-balanced for some integer $n>1$.
\end{Exam}
\noindent{\bf Proof}\ Note that $\{4231,3412\}=\{\overline{1324},\overline{2143}\}$. So it suffices to verify $S_n(1324,2143)$. In addition, it is clear that $S_n(1324,2143)$ is sign-balanced for all integer $1< n\leq4$. Next we observe that $S_{5}(2341,2143)$ is not sign-balanced.

Let $X^{(i)}=\{\pi\in {S_{5}(1324,2143)}:\pi(i)=5\}$. It is clear that $5\pi_1\pi_2\pi_{3}\pi_{4}\in X^{(5)}$ if and only if $\pi_1\pi_2\pi_{3}\pi_{4}\in {S_{4}(1324,2143)}$. Thus $X^{(1)}$ is sign-balanced. Similarly, we deduce that $X^{(2)}$ is sign-balanced. By the computation, it follows that $X^{(3)}=\{12534$, $13542$, $14523$, $23514$, $24531$, $34512$, $41532$, $42513$, $43521\}\cup\{12543$, $14532$, $23541$, $24513$, $34521$, $41523$, $42531$, $43512\}$. Note that the number of even permutations is one more than that of odd permutations in $X^{(3)}$. Similarly, $X^{(4)}=\{12453$, $23451, 31452$, $34251, 41253$, $42351, 43152\}\cup\{12354$, $13452, 32451$, $34152, 41352, 43251\}$ and the number of even permutations is also one more than that of odd permutations in $X^{(4)}$; $X^{(5)}=\{12345$, $23145, 31245$, $32415, 34125$, $42135, 43215\}\cup \{21345$, $23415, 32145$, $34215, 41235$, $42315, 43125\}$ and the number of even permutations is equal to that of odd permutations. Hence, $S_5(1324,2143)$ is not sign-balanced.   $\hfill{\Box}$

%

\vskip0.4cm
\noindent {\bf Acknowledgments}

\noindent The work was supported by the National Natural Science Foundation of China (No. 12061030) and Hainan Provincial Natural Science Foundation of China (No. 122RC652).


{\small


\begin{thebibliography}{99}

\bibitem{A1}
 R.M. Adin, Y. Roichman, Equidistribution and sign-balance on $321$-avoiding permutations, \emph{S$\acute{e}$m. Lothar. Combin.}, {\bf 51} (2004), Artical B51d, 14pp.

\bibitem{A}
M.D. Atkinson, Restricted permutations, \emph{Discrete Math.}, {\bf 195}(1-3) (1999), 27-38.

\bibitem{Baril}
J.-L. Baril, Avoiding patterns in irreducible permutations, \emph{Discrete Math. Theor. Comput. Sci.}, {\bf 17}(3) (2016), 13-30.

\bibitem{BJ}
A. Baxter, A.D. Jaggard, Pattern avoidance by even permutations, \emph{Electron. J. Combin.}, {\bf 18}(2) (2011), 28.


\bibitem{B}
M. B\'ona, \emph{Combinatorics of Permutations}, Second Edition, CRC Press, 2012.


\bibitem{C}
D.I.A. Cohen, \emph{Basic Techniques of Combinatorial Theory}, John Wiley \& Sons, New York, 1978.

\bibitem{EFP}
S.-P. Eu, T.-S. Fu, Y.-J. Pan, A refined sign-balance of simsun permutations, \emph{European J. Combin.}, {\bf 36} (2014), 97-109.

\bibitem{EFP1}
S.-P. Eu, T.-S. Fu, Y.-J. Pan, C.-T. Ting, Sign-balance identities of Adin-Roichman type on $321$-avoiding alternating permutations, \emph{Discrete Math.}, {\bf 312}(15) (2012), 2228-2234.

\bibitem{K}
D.E. Knuth, \emph{The Art of Computer Programming}, Volume $I$: Fundamental Algorithms, Addision-Wesley, 1973.


\bibitem{M1}
T. Mansour, Equidistribution and sign-balance on $132$-avoiding permutations, \emph{S$\acute{e}$m. Lothar. Combin.}, {\bf 51} (2004), Artical B51e, 11pp.

\bibitem{M2}
T. Mansour, Restricted even permutations and Chebyshev polynomials, \emph{Discrete Math.}, {\bf 306}(12) (2006), 1161-1176.

\bibitem{R}
A. Reifegerste, Refined sign-balance on $321$-avoiding permutations, \emph{European J. Combin.}, {\bf 26} (2005), 1009-1018.

\bibitem{SS}
R. Simion, F.W. Schmidt, Restricted permutations, \emph{European J. Combin.}, {\bf 6}(4) (1985), 383-406.


\bibitem{S}
R.P. Stanley, Some remarks on sign-balanced and maj-balanced posets, \emph{Adv. Appl. Math.}, {\bf 34}(4) (2005), 880-902.

\bibitem{V}
V. Vatter, \emph{Permutation Classes}, in: Mikl\'os B\'ona (Ed.), Handbook of Enumerative Combinatorics, CRC Press, 2005.

\end{thebibliography}
\end{document}